\documentclass{article}
\usepackage[T1]{fontenc}
\usepackage{geometry}
\usepackage{amssymb,amsmath,amsthm}
\usepackage{cases}
\usepackage{enumitem}
\usepackage{xcolor}
\usepackage{comment}
\usepackage{graphicx}\graphicspath{{img/}}
\usepackage{textcomp}
\usepackage{mathrsfs}
\usepackage{bbm}
\usepackage{bm}
\usepackage{caption}
\usepackage{stmaryrd}
\usepackage{url}
\usepackage{calc}
\usepackage{ifthen}
\usepackage{psfrag}
\usepackage[active]{srcltx}   
\usepackage{mathtools}
\usepackage{todonotes}
\usepackage[colorlinks=true,urlcolor=magenta!60!black,citecolor=green!50!black,linkcolor=red!90!black,bookmarksnumbered=true,pdfauthor={Jeremie Bettinelli,Dimitri Korkotashvili}]{hyperref}

\theoremstyle{plain}
\newtheorem{thm}{Theorem}
\newtheorem{prop}[thm]{Proposition}
\newtheorem*{claim}{Claim}
\newtheorem{corol}[thm]{Corollary}

\theoremstyle{definition}
\newtheorem{rem}{Remark}

\newcommand{\de}{\mathrel{\mathop:}\hspace*{-.6pt}=}
\newcommand{\un}[1]{\mathbbm{1}_{#1}}
\newcommand\bO{\mathcal{O}}

\newcommand{\vd}{\smash{\vec{\vphantom{i}\smash{\mathsf{d}}}}}
\newcommand{\sew}{\,{\Join}\,}
\newcommand{\LL}{\textup{light-left}}
\newcommand{\LR}{\textup{light-right}}

\newcommand{\ba}{\bm{a}}
\newcommand{\bd}{\bm{d}}
\newcommand{\bl}{\bm{\ell}}
\newcommand{\bgamma}{\boldsymbol{\gamma}}
\newcommand{\bwp}{\boldsymbol{\wp}}

\newcommand{\N}{\mathbb{N}}

\newcommand{\sC}{\mathscr{C}}

\newcommand{\cH}{\mathcal{H}}
\newcommand{\cM}{\mathcal{M}}
\newcommand{\cN}{\mathcal{N}}
\DeclareMathAlphabet\mathbfcal{OMS}{cmsy}{b}{n}
\newcommand{\co}{\scalebox{.6}[.68]{$\mathbfcal{O}$}}

\newcommand{\m}{\mathfrak{m}}

\definecolor{vert}{rgb}{0,0.666,0}
\definecolor{violet}{RGB}{184,113,255}
\definecolor{stepcol}{rgb}{0,.44,.22}

\usepackage{contour}
\contourlength{0.6pt}
\newcommand{\unli}[1]{\underline{\phantom{\smash{\textbf{#1}}}}\llap{\contour{white}{\textbf{#1}}}}

\makeatletter
\renewcommand\paragraph{\@startsection{paragraph}{4}{\z@}%
                                    {3.25ex \@plus1ex \@minus.2ex}%
                                    {-1em}%
                                    {\normalfont\normalsize\bfseries\unli}}
\makeatother

\title{Slit-slide-sew bijections for constellations and quasiconstellations}
\author{J\'er\'emie Bettinelli\thanks{cnrs \& Laboratoire d'Informatique de l'\'Ecole polytechnique; \href{mailto:jeremie.bettinelli@normalesup.org}{\nolinkurl{jeremie.bettinelli@normalesup.org}};\newline \nolinkurl{www.normalesup.org/}\texttildelow\nolinkurl{bettinel}. Partially supported by Grant ANR-20-CE48-0018 \emph{3DMaps}.}%
\and %
Dimitri Korkotashvili\thanks{Laboratoire d'Informatique de l'\'Ecole polytechnique; \href{mailto:dimitri.korkotashvili@polytechnique.edu}{\nolinkurl{dimitri.korkotashvili@polytechnique.edu}}.}}

\includecomment{lvc}
\excludecomment{svc}
\let\oldemph\emph
\renewcommand{\emph}[1]{\textcolor{red!65!black}{\oldemph{#1}}}

\begin{document}
\maketitle

\begin{abstract}
We extend so-called slit-slide-sew bijections to constellations and quasiconstellations, which allow to recover the counting formula for constellations or quasiconstellations with a given face degree distribution.

More precisely, we present an involution on the set of hypermaps given with an orientation, one distinguished corner, and one distinguished edge leading away from the corner while oriented in the given orientation. This involution reverts the orientation, exchanges the distinguished corner with the distinguished edge in some sense, slightly modifying the degrees of the incident faces in passing, while keeping all the other faces intact.

The construction consists in building a canonical path from the distinguished elements, slitting the map along it, and sewing back after sliding by one unit along the path. The involution specializes into a bijection interpreting combinatorial identities linking the numbers of constellations or quasiconstellations with a given face degree distribution, where the degree distributions differ by one $+1$ and one $-1$.

Our bijections yield a ``degree-by-degree, face-by-face'' growth algorithm that samples a hypermap uniformly distributed among constellations or quasiconstellations with prescribed face degrees. More precisely, it samples at each step uniform constellations or quasiconstellations, whose face degree distributions slightly evolve to the desired distribution.
\end{abstract}


\section{Introduction}

In the present work, which is the extended version of~\cite{BeKo24short}, we pursue the investigation of so-called \emph{slit-slide-sew} bijections, introduced in~\cite{Bet14} on forests and plane quadrangulations, and then further developed in \cite{Bet19short,Bet20} on plane bipartite and quasibipartite maps. Here, we focus on a generalization of the latter, called \emph{constellations} and \emph{quasiconstellations}.

Constellations have been the focus on several studies. They bear connections with factorizations of permutations and branched coverings of the sphere (see e.g.\ \cite{BMSc00}). They are also deeply related to so-called Hurwitz numbers~\cite{Hur91,Eze78,DuPoSc14}. See for instance~\cite{BoChChGF24} and the references therein for a broad overview of these developments.

\paragraph{Hypermaps.}
Recall that a \emph{plane map} is an embedding of a finite connected graph (possibly with multiple edges and loops) into the sphere, considered up to orientation-preserving homeomorphisms. Now fix an integer $p\ge 2$. A (plane) \emph{$p$-hypermap} is a plane map whose faces are shaded either \emph{dark} or \emph{light} in such a way that 
\begin{itemize}
	\item adjacent faces do not have the same shade (one is dark, the other light);
	\item each dark face has degree~$p$.
\end{itemize}
These actually generalize maps, which correspond to $2$-hypermaps. In the terminology of hypermaps, light faces generalize faces and might be called \emph{hyperfaces}, whereas dark faces generalize edges and are called \emph{hyperedges}. We do not use this terminology here. Hypermaps may alternatively be defined in terms of permutations: see~\cite[Section~2]{BMSc00}.

\begin{svc}
\begin{wrapfigure}[17]{r}{0.49\textwidth}
		\psfrag{f}[][][.9]{$f_1$}
		\psfrag{a}[][][.9]{$f_4$}
		\psfrag{z}[][][.9]{$f_8$}
		\psfrag{e}[][][.9]{$f_3$}
		\psfrag{r}[][][.9]{$f_5$}
		\psfrag{t}[][][.9]{$f_7$}
		\psfrag{y}[][][.9]{$f_6$}
		\psfrag{u}[][][.9]{$f_9$}
		\psfrag{i}[][][.9]{$f_2$}
	\centering\includegraphics[width=.95\linewidth]{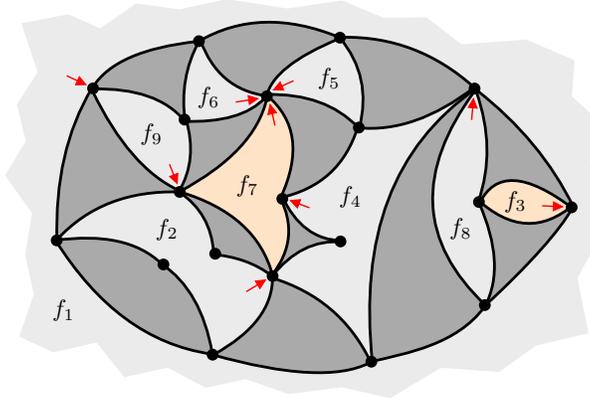}
	\vspace{-5mm}
	\caption{A quasi-$3$-constellation of type $(9,6,2,6,3,3,4,3,3)$. The two flawed faces are~$f_3$ and~$f_7$; throughout the paper, we highlight flawed faces in orange. Every light face has a marked corner, always represented by a red arrowhead.}
	\label{phypermap}
\end{wrapfigure}
\end{svc}

A (plane) \emph{$p$-constellation} is a $p$-hypermap such that the degrees of its light faces are all multiples of~$p$. In a \emph{$p$-hypermap}, a light face whose degree is not a multiple of~$p$ will be called a \emph{flawed face}. A $p$-constellation is thus a $p$-hypermap without flawed faces. A \emph{quasi-$p$-constellation} is a $p$-hypermap with exactly two flawed faces. Note that, in a $p$-hypermap, the sum of the degrees of the light faces is necessarily a multiple of~$p$, since it is equal to the sum of the degrees of the dark faces, which are all~$p$. As a result, a $p$-hypermap cannot have a single flawed face and, in a quasi-$p$-constellation, the two flawed faces have, modulo~$p$, degrees~$+k$ and~$-k$ for some $0< k <p$.

\begin{figure}[ht]
		\psfrag{f}[][][.9]{$f_1$}
		\psfrag{a}[][][.9]{$f_4$}
		\psfrag{z}[][][.9]{$f_8$}
		\psfrag{e}[][][.9]{$f_3$}
		\psfrag{r}[][][.9]{$f_5$}
		\psfrag{t}[][][.9]{$f_7$}
		\psfrag{y}[][][.9]{$f_6$}
		\psfrag{u}[][][.9]{$f_9$}
		\psfrag{i}[][][.9]{$f_2$}
	\centering\includegraphics[width=8cm]{phypermap}
	\caption{A quasi-$3$-constellation of type $(9,6,2,6,3,3,4,3,3)$. The two flawed faces are~$f_3$ and~$f_7$; throughout the paper, we highlight flawed faces in orange. Every light face has a marked corner, always represented by a red arrowhead.}
	\label{phypermap}
\end{figure}

\begin{figure}[ht]
		\psfrag{f}[][][.9]{$f_1$}
		\psfrag{a}[][][.9]{$f_5$}
		\psfrag{z}[][][.9]{$f_8$}
		\psfrag{e}[][][.9]{$f_3$}
		\psfrag{r}[][][.9]{$f_4$}
		\psfrag{t}[][][.9]{$f_6$}
		\psfrag{y}[][][.9]{$f_7$}
		\psfrag{u}[][][.9]{$f_9$}
		\psfrag{i}[][][.9]{$f_2$}
	\centering\includegraphics[width=.9\linewidth]{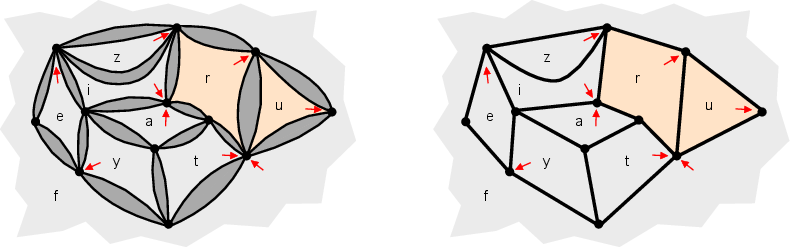}
	\caption{A quasi-$2$-constellation and the corresponding quasibipartite map, obtained by squeezing all the dark faces.}
	\label{squeeze}
\end{figure}

In the particular case $p=2$, squeezing each dark face (which has degree~$2$ by definition) into a single edge yields a simple one-to-one correspondence between $2$-constellations and bipartite maps, as well as between quasi-$2$-constellations and quasibipartite maps (maps with exactly two odd degree faces), as shown in Figure~\ref{squeeze}.

Another particular class of interest enters our framework: the class of Eulerian triangulations, that is, $3$-constellations whose light faces are all of degree~$3$.

\paragraph{Enumeration.}
For an $r$-tuple $\ba=( a_1,\ldots, a_r)$ of positive integers, let us denote by $C(\ba)$ the number of $p$-hypermaps with exactly~$r$ light faces, numbered~$f_1$, \ldots, $f_r$ and of respective degrees~$a_1$, \ldots, $ a_r$, each bearing a marked corner\footnote{Recall that a \emph{corner} is an angular sector delimited by two consecutive edges around a vertex.}. The $r$-tuple~$\ba$ will be called the \emph{type} of such $p$-hypermaps. See Figure~\ref{phypermap}. By elementary considerations and Euler's characteristic formula, the integers
\begin{equation*}
E(\ba)\de \sum_{i=1}^r a_i\,,\qquad D(\ba)\de\frac{E(\ba)}p \,,\qquad \text{ and }\qquad V(\ba)\de E(\ba)- D(\ba) - r + 2
\end{equation*}
are respectively the numbers of edges, dark faces, and vertices of $p$-hypermaps of type~$\ba$. Generalizing Tutte's so-called \emph{formula of slicings}~\cite{Tut62sli}, it has been computed that, when at most two~$a_i$'s are not in~$p\N$, that is, for $p$-constellations~\cite{BMSc00} or quasi-$p$-constellations~\cite{CoFu14}, it holds that

\begin{align}
C(\ba)=c_{\ba}\,\frac{\big(E(\ba)-D(\ba)-1\big)!}{V(\ba)\,!}\prod_{i=1}^r\alpha( a_i),\qquad
	\text{ where }	&\quad\alpha(x)\de\frac{x!}{\big\lfloor x/p\big\rfloor!\,\big(x-\lfloor x/p\rfloor-1\big)!}\label{slicings}\\\notag
	\text{ and}	&\quad
c_{\ba}=\begin{cases}
	1 	&\text{ if $p$ divides every~$a_i$}\\
	p-1 	&\text{ otherwise}
\end{cases}.
\end{align}

\paragraph{Combinatorial identities.}
In the present work, we give a bijective interpretation for the following combinatorial identity, which transfers one degree from one face to another.

\begin{prop}[Transferring one degree from~$f_1$ to~$f_2$]\label{propmp1}
Let $\ba=(a_1,\ldots, a_r)$ be an $r$-tuple of positive integers such that $a_1\ge 2$, and with coordinates congruent modulo~$p$ to
\begin{enumerate}[label=(\textit{\roman*})]
	\item\label{congrui} either $(k,-k,0,\dots,0)$ for some $k\in\{0,\dots,p-1\}$,
	\item\label{congruii} or $(1,0,\dots, 0,-1,0,\dots,0)$, with the~$- 1$ in any position from~$3$ to~$r$.
\end{enumerate}
Let also $\tilde \ba=(\tilde a_1,\ldots,\tilde a_r)\de( a_1-1, a_2+1, a_3,\ldots, a_r)$. Then the following identity holds:
\begin{align}\label{eqpm}
\big(a_1 - \lceil a_1/p\rceil \big)\, \big(a_2+1\big)\, C(\ba)&=\big(\tilde a_1+1\big)\,\big(\tilde a_2 - \lceil \tilde a_2/p\rceil \big)\, C(\tilde \ba)\,.
\end{align}
\end{prop}

To obtain~\eqref{eqpm} from~\eqref{slicings}, one might first observe that, for any $x\in \N$,
\[
\frac{\alpha(x)}{\alpha(x-1)}=d_x\,\frac{x}{x- \lceil x / p\rceil}\qquad\text{ where }\quad
d_x=\begin{cases}
	p-1 & \text{ if } p \mid x\\
	1 & \text{ if } p \nmid x 
\end{cases},
\]
and then that, in both cases~\ref{congrui} and~\ref{congruii}, $c_{\ba} d_{a_1} = c_{\tilde\ba} d_{\tilde a_2} = p-1$.

\begin{rem}
In fact, Proposition~\ref{propmp1} also holds when the coordinates of~$\ba$ are congruent modulo~$p$ to
\begin{enumerate}[label=(\textit{\roman*}),start=3]
	\item\label{congruiii} $(0,-1,0,\dots, 0,1,0,\dots,0)$, with the~$+1$ in any position from~$3$ to~$r$,
\end{enumerate}
observing that, in this case, $c_{\ba} d_{a_1} = c_{\tilde\ba} d_{\tilde a_2}=(p-1)^2$. For some technical reason, however, this particular case falls out of the reach of the present paper.
\end{rem}

We furthermore treat the case of a degree $1$-face, which may easily be obtained as above.

\begin{prop}[Transferring the degree of a degree~$1$-face~$f_1$ to~$f_2$]\label{propface}
Let $\ba=(1,a_2,\ldots, a_r)$ and $\tilde \ba=(\tilde a_2,\ldots,\tilde a_r)\de(a_2+1, a_3,\ldots, a_r)$ be respectively an $r$-tuple and an $r-1$-tuple of positive integers, both having at most two coordinates not lying in~$p\N$. Then the following identity holds:
\begin{align}\label{eqpm0}
\big(a_2+1\big)\, C(\ba)&=V(\tilde\ba)\,\big(\tilde a_2 - \lceil \tilde a_2/p\rceil \big)\, C(\tilde \ba)\,.
\end{align}
\end{prop}

\begin{figure}[ht!]
		\psfrag{f}[][][.9]{$f_1$}
	\centering\includegraphics[width=.99\linewidth]{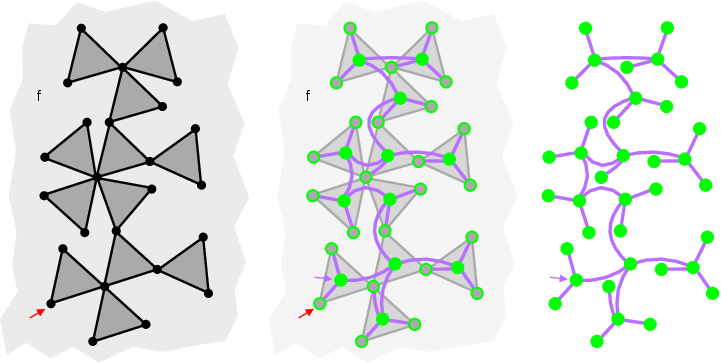}
	\caption{Bijection between $p$-constellations of type $(pn)$ and $p$-ary trees with~$n$ nodes. Each dark face corresponds to a $p$-ary node and each vertex to a leaf. We follow the contour of the light face, starting from its marked corner, in a given orientation, say with the light face on the left. Each time we encounter an unvisited dark face, say~$D$, we consider the~$3$ incident vertices, say~$v_1$, $v_2$, $v_3$. For $1\le i \le 3$, we link the node corresponding to~$D$ to the node corresponding to the left-most unvisited dark face incident to~$v_i$ or, by default, to the leaf corresponding to~$v_i$. We root the tree at the corner following the marked corner of~$f_1$.}
	\label{bijpary}
\end{figure}

\paragraph{Recovering the counting formula~\eqref{slicings}.}
It is easy to see that the number of $p$-constellations with exactly one light face of degree~$pn$ is equal to the number of $p$-ary trees with~$n$ nodes: see Figure~\ref{bijpary}. It thus holds (see for instance~\cite[Section~2.7.7.]{GoJa83}) that
\[
C\big((pn)\big)=\frac{(pn)!}{n!\, \big((p-1)n+1\big)!}.
\]
Using this formula as initial condition, Propositions~\ref{propmp1} and~\ref{propface} provide yet another proof of~\eqref{slicings}.

\paragraph{Methodology.}
In order to bijectively interpret the combinatorial identities~\eqref{eqpm} and~\eqref{eqpm0}, the idea is to distinguish elements, such as edges, vertices, faces, corners, etc., in such a way that each side of an equation of interest counts maps given with such distinguished elements. For instance, in the left-hand side of~\eqref{eqpm}, the term $a_2+1$ counts the number of corners\footnote{The $+1$ comes from the presence of the marked corner in the face. See the first paragraph of Section~\ref{secprel} for the convention we use here.} in the face~$f_2$ of a $p$-hypermap of type~$\ba$. As a result, $\big(a_2+1\big)\, C(\ba)$ is the number of $p$-hypermaps of type~$\ba$ with a distinguished corner in the face~$f_2$. The other prefactor in the left-hand side of~\eqref{eqpm} is more delicate to interpret; this will be done in due course. Remark that we will always use the word ``distinguished'' to designate these extra elements, keeping the word ``marked'' only for the marked corners, which we see as inherent to the considered hypermaps.

Once both sides of the considered equation are properly interpreted as cardinalities of sets of maps with distinguished elements, we bijectively go from one set to the other as follows. Using the distinguished elements, we construct a directed path in the map, called \emph{sliding path}. We then slit the map along this sliding path and sew back together the sides of the slit after sliding by one unit, in the sense that the left side of the $i$-th edge is sewn back on the right side of the $i\pm 1$-th edge (the $\pm1$ being the same for all edges and determined by some rule). This mildly modifies the map along the path but does not affect its faces, except the two that are around the extremities of the sliding path. In the process, new distinguished elements naturally appear in the resulting map; these allow us to recover the sliding path in order to slide back.

\paragraph{Growth algorithm.}
An important feature of our bijections is that they allow to ``grow'' uniform $p$-constellations or quasi-$p$-constellations up to any given type $\ba=(a_1,\ldots, a_r)$. Without loss of generality, let us assume that~$p$ divides every~$a_i$ for $1\le i \le r-2$. In order to address this question, we first denote by $n\de D(\ba)$ the number of dark faces (so that $np=a_1+a_2+\dots+a_r$) and define the finite sequence of tuples: 
\begin{equation}\label{eqseq}
\begin{aligned}
(np),\\
(1,np-1),\dots,(a_1,np-a_1),\\
(a_1,1,np-a_1-1),\dots,(a_1,a_2,np-a_1-a_2),\\
(a_1,a_2,1,np-a_1-a_2-1),\dots,(a_1,a_2,a_3,np-a_1-a_2-a_3)\\
\vdots\hspace{30mm}\\
(a_1,a_2,\dots,a_{r-2},1,a_r+a_{r-1}-1),\dots,(a_1,a_2,\dots,a_r),
\end{aligned}
\end{equation}
which we denote by~$\ba^{0}$, \ldots, $\ba^{np-a_r}=\ba$ for simplicity. The important point is that each tuple of the sequence has at most two coordinates not lying in~$p\N$ and, for each $1 \leq j \leq np-a_r$, the tuple~$\ba^{j}$ differs from~$\ba^{j-1}$ by~$-1$ in one coordinate and
\begin{enumerate}[label=(\textit{\alph*})]
    \item\label{samplea} either~$+1$ in another coordinate,
    \item\label{sampleb} or by a new coordinate equal to~$1$,
\end{enumerate}
while all unchanged coordinates being multiples of~$p$. In particular, up to an obvious coordinate permutation, this change in coordinates is covered by Proposition~\ref{propmp1}.\ref{congrui} or by Proposition~\ref{propface}. We chose this particular sequence but, in fact, any other sequence with the above properties also works. Note also that the number of dark faces corresponding to any of the tuples in the sequence is equal to~$n$.

\medskip
Using the bijection of Figure~\ref{bijpary}, sampling a uniform $p$-constellation~$\m^0$ of type $\ba^0=(np)$ amounts to sampling a uniform $p$-ary tree with~$n$ nodes, which can be done by sampling an appropriate {\L}ukasiewicz path (see~\cite{AlReSc97} for instance).

\begin{figure}[ht!]
		\setlength{\fboxsep}{0pt}
		\definecolor{facegreen}{RGB}{142,255,142}
		\definecolor{faceyellow}{RGB}{255,255,199}
		\definecolor{facepink}{RGB}{255,199,255}
		\definecolor{faceblue}{RGB}{142,255,255}
		\psfrag{a}[B][t]{$(\colorbox{facepink}{$21$})$}
		\psfrag{b}[B][t]{$(\colorbox{faceyellow}{$1$},\colorbox{facepink}{$20$})$}
		\psfrag{c}[B][t]{$(\colorbox{faceyellow}{$6$},\colorbox{facepink}{$15$})$}
		\psfrag{d}[B][t]{$(\colorbox{faceyellow}{$6$},\colorbox{facegreen}{$3$},\colorbox{facepink}{$12$})$}
		\psfrag{e}[B][t]{$(\colorbox{faceyellow}{$6$},\colorbox{facegreen}{$3$},\colorbox{faceblue}{$1$},\colorbox{facepink}{$11$})$}
		\psfrag{x}[B][t]{$\dots$}
		\vspace{3mm}
	\centering\includegraphics[width=.95\linewidth]{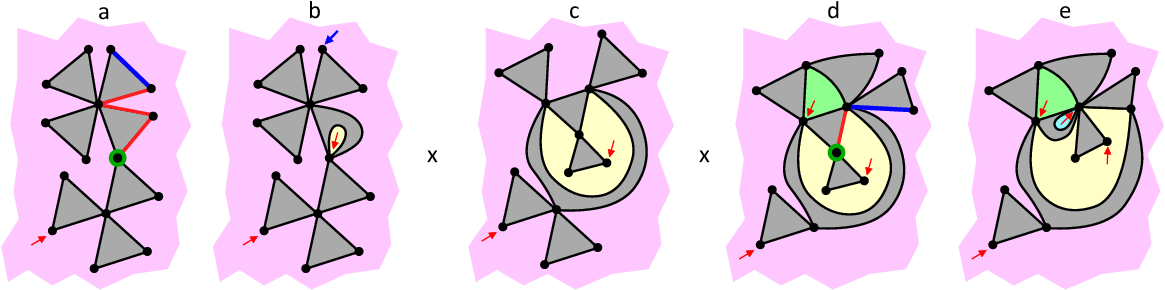}
	\caption{Growing a $3$-constellation of type $(6,3,1,11)$. \textbf{Step 0.} We start with a random $3$-constellation of type $(21)$, thus having a unique (purple) light face, which we think of an external face. \textbf{Step 1.} Randomly selecting a (green) vertex and a (blue) edge, our bijections allow to add a new degree-$1$ (yellow) light face at the cost of removing one degree from the external face. The resulting quasi-$3$-constellation comes with a distinguished (blue) corner in the external face, which we forget for the next step. \textbf{Step 2--6.} We repeat the process of randomly transferring degrees from the external light face to the yellow face, thus obtaining a random $3$-constellation of type $(6,15)$. \textbf{Step 7--9.} We add a new (green) degree-$1$ light face then ``inflate'' it to obtain a random $3$-constellation of type $(6,3,12)$. \textbf{Step 10.} Finally, we add a new (blue) degree-$1$ light face to obtain a random quasi-$3$-constellation of type $(6,3,1,11)$. In this final step, the selected (green) vertex and (blue) edge used for applying the bijection are also represented.} 
	\label{figrow}
\end{figure}

Then, as we will explain during Section~\ref{secsample} in more details, our bijections allow to sample a sequence of $p$-hypermaps $\m^0$, \dots, $\m^{np-a_r}$ such that, for every $1\le j\le np-a_r$, the $p$-hypermap~$\m^j$ is uniformly distributed among $p$-hypermaps of type~$\ba^j$. See Figure~\ref{figrow}.

Note that, in the sequence~\ref{eqseq}, the $(i+1)$-th line adds a new face of degree growing from~$1$ to its desired degree~$a_i$, while the last face sees its degree decrease by~$a_i$. Seeing this last face as an \emph{external face}, we obtain a sequence of ``growing'' uniform $p$-hypermaps where the faces are added one by one, each face being grown by adding one degree at each step. In other words, we start from a uniform tree made of~$n$ dark degree $p$-faces. Then we grow a first face of degree~$1$, then~$2$, and so on until it has degree~$a_1$, then we grow a second face of degree~$1$, then~$2$, and so on until it has degree~$a_2$, and we keep growing faces until we reach the desired type~$\ba$.

\paragraph{Organization of the paper.}
The remainder of the document is structured in the following manner. We start by giving in Section~\ref{secprel} the definitions and conventions we use, as well as a combinatorial interpretation of the prefactor $\big(a - \lceil a/p\rceil \big)$ appearing in the identities~\eqref{eqpm} and~\eqref{eqpm0}. We then present in Section~\ref{secbij} our bijective interpretation of these identities through a more general involution on the set of maps given with an orientation, a distinguished corner, and a distinguished edge satisfying an extra constraint. Finally, we explain in Section~\ref{secsample} how to sample a uniform $p$-hypermap of a given type using our bijections.

\section{Preliminaries}\label{secprel}

\begin{svc}
\begin{wrapfigure}[5]{r}{0.55\textwidth}
		\psfrag{c}[][]{\textcolor{blue}{$c$}}
		\psfrag{o}[][]{or}
	\vspace{-8mm}
	\centering\includegraphics[width=6cm]{discormarksv}
	\caption{Distinguishing a corner around the marked corner.}
	\vspace{-5mm}
	\label{discormark}
\end{wrapfigure}
\end{svc}
	
\paragraph{Distinguishing a corner.}
Following previous works on slit-slide-sew bijections, we use the convention, depicted in Figure~\ref{discormark}, that the marked corner of a face creates two possible corners to distinguish. One might think of the marked corner as a dangling half-edge, with one corner on each side. As a result, a face of degree~$a$ bearing its marked corner has $a+1$ possible corners to distinguish.

\begin{figure}[ht]
		\psfrag{c}[][]{\textcolor{blue}{$c$}}
		\psfrag{o}[][]{or}
	\centering\includegraphics[width=10cm]{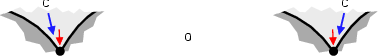}
	
	\caption{Distinguishing a corner around the marked corner.}
	\label{discormark}
\end{figure}

\paragraph{Edge orientation.}
As is customary, we will orient the edges of the hypermaps we consider, in such a way that light faces always lie to the same side of the oriented edges (and thus dark faces always lie to the same other side). These orientations will be called the \emph{light-left orientation} when the light faces\footnote{Recall that the light faces are the main objects of focus.} all lie to the left, and the \emph{light-right orientation} when the light faces all lie to the right. In other words, in the light-right orientation, the edges are oriented clockwise around light faces and counterclockwise around dark faces. See Figure~\ref{orien}. We will need to use both orientations in the present paper. We will always clearly mention which orientation we use whenever it matters. Without specific mention, both orientations can be used. \textbf{Once one of the two possible orientations is fixed}, we will use the following conventions.

\begin{svc}
\begin{wrapfigure}[11]{r}{0.45\textwidth}
		\psfrag{a}[][][.9]{$f_5$}
		\psfrag{r}[][][.9]{$f_4$}
		\psfrag{t}[][][.9]{$f_6$}
		\psfrag{y}[][][.9]{$f_7$}
		\psfrag{u}[][][.9]{$f_9$}		
		\psfrag{e}[][]{\textcolor{vert}{$e$}}
		\psfrag{-}[][]{\textcolor{vert}{$e^-$}}
		\psfrag{+}[][]{\textcolor{vert}{$e^+$}}
		\psfrag{b}[][]{\textcolor{vert}{$c_e$}}
		\psfrag{c}[][]{\textcolor{blue}{$c$}}
		\psfrag{d}[][]{\textcolor{blue}{$c^+$}}
	\vspace{-3mm}
	\centering\includegraphics[width=7cm]{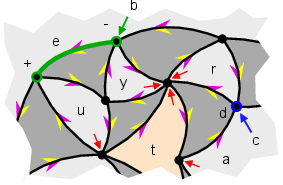}
	\vspace{-5mm}
	\caption{Edge orientation and related definitions. Here, the light-right orientation is depicted.}
	\label{orien}
\end{wrapfigure}
\end{svc}

Given an edge~$e$, we will respectively denote by~$e^{-}$ and~$e^{+}$ the origin and end of the edge~$e$, oriented as convened. The corner \emph{preceding}~$e$ is defined as the corner~$c_e$ delimited by~$e$ and the edge that precedes~$e$ in the contour of the incident light face, in the convened orientation. Similarly, we denote by~$c^+$ the vertex incident to a corner~$c$.

\begin{figure}[ht!]
		\psfrag{a}[][][.9]{$f_5$}
		\psfrag{r}[][][.9]{$f_4$}
		\psfrag{t}[][][.9]{$f_6$}
		\psfrag{y}[][][.9]{$f_7$}
		\psfrag{u}[][][.9]{$f_9$}		
		\psfrag{e}[][]{\textcolor{vert}{$e$}}
		\psfrag{-}[][]{\textcolor{vert}{$e^-$}}
		\psfrag{+}[][]{\textcolor{vert}{$e^+$}}
		\psfrag{b}[][]{\textcolor{vert}{$c_e$}}
		\psfrag{c}[][]{\textcolor{blue}{$c$}}
		\psfrag{d}[][]{\textcolor{blue}{$c^+$}}
	\centering\includegraphics[width=7.4cm]{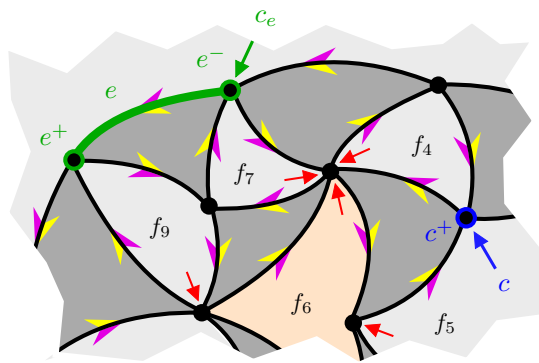}
	\caption{Edge orientation and related definitions. Here, the light-right orientation is depicted.}
	\label{orien}
\end{figure}

\paragraph{Paths.}
A \emph{path} from a vertex~$v$ to a vertex~$v'$ is a finite sequence $\bwp=(e_1,e_2,\dots,e_k)$ of edges such that $e_1^-=v$, for $1\le i \le k-1$, $e_i^+=e_{i+1}^-$, and $e_k^+=v'$. Its \emph{length} is the integer~$k$, which we denote by $[\bwp]\de k$. A path is called \emph{simple} if the vertices it visits are all distinct.

Beware that a path is only made of edges oriented in the convened orientation. In other words, edges cannot be used ``backward.'' In particular, this means that all the faces lying to the left of a path are of the same shade (either all light or all dark), whereas all the faces lying to its right are of the other shade. The side of the path where the faces are all light will be called its \emph{light side}, whereas the other side will be called its \emph{dark side}.

\paragraph{Directed metric and geodesics.}
We will use the directed metric associated with the convened orientation: given two vertices~$v$, $v'$ in a $p$-hypermap, we denote by $\vd(v, v')$ the 
smallest~$k$ for which there exists a path from~$v$ to~$v'$ of length~$k$. (We put an arrow on top in the notation to keep in mind that this is only a directed metric.) A \emph{geodesic} from~$v$ to~$v'$ is such a path. 

\begin{svc}
\begin{wrapfigure}[11]{r}{0.35\textwidth}
		\psfrag{k}[][]{\textcolor{blue}{$e_k$}}
		\psfrag{+}[l][c]{\textcolor{red}{$e_{k+1}$}}
	\vspace{-5mm}
	\centering\includegraphics[width=4cm]{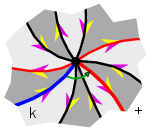}
	\vspace{-5mm}
	\caption{Definition of the lightest geodesic. The edges going closer to~$v'$ are in red.}
	\label{figgeod}
\end{wrapfigure}
\end{svc}

There are generally several geodesics from a given vertex~$v$ to a target vertex~$v'$. Among all of these, one will be of particular interest in this work: the \emph{lightest geodesic}, constructed as follows. It is only well defined from a starting edge or corner~$e_0$ such that $e_0^+=v$. (The starting element~$e_0$ does not belong to the path.) Then, provided~$e_0$, $e_1$, \dots, $e_j$ have already been constructed and the path is not complete (that is, $e_j^+\ne v'$), we set the subsequent edge~$e_{j+1}$ as the one, among the edges~$e$ such that $e^-=e_j^+$ and $\vd(e^+,v')=\vd(e_{j}^+,v')-1$, that comes first while turning around~$e_j^+$ in the direction where we successively see ``incoming edge, light face, outgoing edge, dark face''. See Figure~\ref{figgeod}. In other words, the lightest geodesic is the leftmost geodesic if the convened orientation is the light-left orientation and the rightmost geodesic if the convened orientation is the light-right orientation.

\begin{figure}[ht!]
		\psfrag{k}[][]{\textcolor{blue}{$e_k$}}
		\psfrag{+}[l][c]{\textcolor{red}{$e_{k+1}$}}
	\centering\includegraphics[width=4cm]{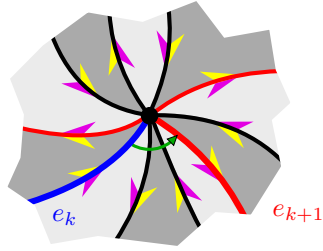}
	\caption{Definition of the lightest geodesic. The three red edges are the ones going closer to the target vertex~$v'$. Turning around~$e_k^+$ in the direction of the light face, we take the first red edge (drawn with a thicker line).}
	\label{figgeod}
\end{figure}

\begin{figure}[ht!]
		\psfrag{0}[B][B][.7]{\textcolor{white}{$0$}}
		\psfrag{1}[B][B][.7]{\textcolor{white}{$1$}}
		\psfrag{2}[B][B][.7]{\textcolor{white}{$2$}}
		\psfrag{3}[B][B][.7]{\textcolor{white}{$3$}}
		\psfrag{4}[B][B][.7]{\textcolor{white}{$4$}}
		\psfrag{5}[B][B][.7]{\textcolor{white}{$5$}}
		\psfrag{6}[B][B][.7]{\textcolor{white}{$6$}}
		\psfrag{v}[][][.8]{\textcolor{vert}{$v$}}
		\psfrag{f}[][][.8]{$f$}
	\centering\includegraphics[width=8cm]{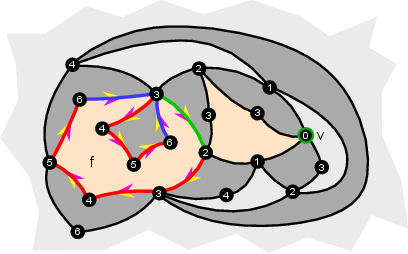}
	\caption{The different types of edges incident to a flawed face in a quasi-$4$-constellation. We use the light-right orientation; the distances to~$v$ are written in the vertices. Around~$f$, the $\big(10 - \lceil 10/4\rceil \big)=7$ red edges are leaving~$v$; the $\lfloor 10/4\rfloor=2$ blue edges are approaching~$v$; the green edge is irregular with respect to~$v$.}
	\label{edgetyp}
\end{figure}

\paragraph{Edge types.}
Given a fixed vertex~$v$ in a $p$-hypermap, we may differentiate three types of edges: an edge~$e$ is said to be 
\begin{itemize}
	\item \emph{leaving~$v$} if $\vd(v, e^+)=\vd(v, e^-)+1$;
	\item \emph{approaching~$v$} if $\vd(v, e^+)=\vd(v, e^-)+1-p$;
	\item \emph{irregular with respect to~$v$} if $\vd(v, e^+)-\vd(v, e^-) \not\equiv 1 \bmod p$.
\end{itemize}
Observe that $1-p \le \vd(v, e^+)-\vd(v, e^-) \le 1$ since there is always a path of length~$1$, namely the path consisting of the single edge~$e$, as well as a path from~$e^+$ to~$e^-$ of length $p-1$, made of all the other edges incident to the dark face incident to~$e$. As a result, if~$e$ is irregular with respect to~$v$, then it holds that $\vd(v, e^+)-\vd(v, e^-) \in \{2-p,3-p,\dots,0\}$. See Figure~\ref{edgetyp}.

The following proposition bounds the number of irregular edges in a $p$-hypermap in terms of the number of flawed faces.

\begin{prop}\label{propirreg}
We consider a vertex~$v$ and a light face~$f$ in a $p$-hypermap having exactly~$k$ flawed faces. Then, the number of edges incident to~$f$ that are irregular with respect to~$v$ is
\begin{itemize}
	\item at most~$k-1$ if $f$ is a flawed face;
	\item at most~$k$ otherwise.
\end{itemize}
\end{prop}

In the particular cases that are of interest to us, namely for a given face in a $p$-constellation or a given \textbf{flawed} face in a quasi-$p$-constellation, it gives the number of each type among the incident edges. This provides an interpretation to the prefactor $\big(a - \lceil a/p\rceil \big)$ appearing in~\eqref{eqpm} and~\eqref{eqpm0}.

\begin{corol}\label{propq}
We consider a vertex~$v$ and a light face~$f$ of degree~$a$ in a $p$-hypermap.
\begin{enumerate}[label=(\textit{\arabic*})]
\item\label{propqi} If the $p$-hypermap is a $p$-constellation then, among the~$a$ edges incident to~$f$, 
\begin{itemize}
	\item $a - a/p$ are leaving~$v$;
	\item $a/p$ are approaching~$v$;
	\item none are irregular with respect to~$v$.
\end{itemize}

\item\label{propqii} If the $p$-hypermap is a quasi-$p$-constellation and~$f$ a flawed face then, among the~$a$ edges incident to~$f$, 
\begin{itemize}
	\item $\big(a - \lceil a/p\rceil \big)$ are leaving~$v$;
	\item $\lfloor a/p\rfloor$ are approaching~$v$;
	\item one is irregular with respect to~$v$.
\end{itemize}
\end{enumerate}
\end{corol}

\begin{proof}
First of all, observe that, for either~\ref{propqi} or~\ref{propqii}, the third subitem easily implies the first and second subitems. Indeed, by definition and from the observation above, following the contour of~$f$, the (oriented) distance variations to~$v$ are
\begin{itemize}
	\item $+1$ for leaving edges;
	\item $+1-p$ for approaching edges;
	\item $-k$ for some $k\in \{0,1,\dots,p-2\}$ for the irregular edge, in the quasi-$p$-constellation case.
\end{itemize}
These differences must sum up to~$0$, which gives the desired result (and shows in passing that, in fact, $k=a-p\lfloor a/p\rfloor-1$).

\smallskip

Now, note that, in the quasi-$p$-constellation case, there must be at least one irregular edge since, for all other edges, the distance variations to~$v$ are congruent to $1 \bmod p$ and sum up to~$0$. The third subitem comes thus from Proposition~\ref{propirreg}, in the particular case $k=0$ for~\ref{propqi} and $k=2$ for~\ref{propqii}.
\end{proof}

We proceed to the proof of the more general proposition.

\begin{proof}[Proof of Proposition~\ref{propirreg}]
Let us denote by~$e_1$, \dots, $e_m$ the edges incident to~$f$ that are irregular with respect to~$v$. Then, for each $1\le i\le m$, let us consider two geodesics~$\bgamma_i^-$ and~$\bgamma_i^+$ from~$v$ towards~$e_i^-$ and~$e_i^+$. Note that, since~$e_i$ is not leaving~$v$, the geodesics~$\bgamma_i^-$ and~$\bgamma_i^+$ have to split before reaching~$e_i^-$. For any two geodesics, if they meet again after they split, the separate parts inbetween the spliting point and the meeting point have, by definition of geodesicvs, the same length. Therefore, up to changing the geodesics, one may furthermore assume that any two of the considered geodesics never meet again after the time they split. These~$2m$ geodesics make up a directed tree whose leaves are the extremities of~$e_1$, \dots, $e_m$. See Figure~\ref{pfirreg}.

\begin{figure}[ht]
		\psfrag{f}[][]{$f$}
		\psfrag{v}[][][.8]{$v$}
		\psfrag{1}[][]{\textcolor{vert}{$e_1$}}
		\psfrag{2}[][]{\textcolor{vert}{$e_2$}}
		\psfrag{3}[][]{\textcolor{vert}{$e_3$}}
		\psfrag{-}[][]{\textcolor{red}{$\bwp_2^-$}}
		\psfrag{+}[c][t]{\textcolor{red}{$\bwp_2^+$}}
	\centering\includegraphics[width=.85\linewidth]{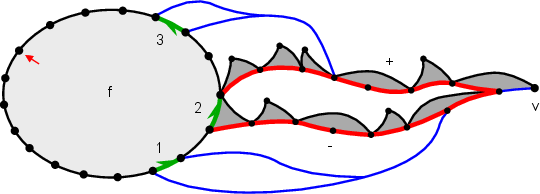}
	\caption{The face~$f$ and the irregular edges~$e_1$, $e_2$, $e_3$. The geodesics are in blue, except the parts~$\bwp_2^-$ and~$\bwp_2^+$, which are highlighted in red. Only the dark faces that are to the dark side of~$\bwp_2^-$ and of~$\bwp_2^+$ are represented.}
	\label{pfirreg}
\end{figure}

Now, for each~$i$, let us consider the parts~$\bwp_i^-$ and~$\bwp_i^+$ of~$\bgamma_i^-$ and~$\bgamma_i^+$ after the vertex where they split up. The concatenation of~$\bwp_i^-$, $e_i$, and then the reverse of~$\bwp_i^+$ provides a cycle of edges (not all oriented in the same direction) in the $p$-hypermap. By the Jordan curve theorem, this cycle separates the sphere into two connected components: let us denote by~$\sC_i$ the one that does not contain the face~$f$. Note that, by construction, $\sC_1$, \dots, $\sC_m$, $f$, are pairwise disjoint. Consequently, the proof of the claim will be complete if we show that each~$\sC_i$ contains at least one flawed face, since this yields that $m+\un{\{f\text{ is flawed}\}} \le k$.

In order to see this, let us consider the $p$-hypermap constructed as shown in Figure~\ref{pfirreg2}. We start from~$\sC_i$ and we add one dark face to the dark side of each of the~$[\bwp_i^+]$ edges on the part of the boundary of~$\sC_i$ that corresponds to~$\bwp_i^+$. Then we add a light face (of degree $[\bwp_i^-]+1+(p-1)[\bwp_i^+]$) on the boundary of the resulting object. (Recall that~$\bwp_i^+$ is taken backward in the cycle and thus has its incident dark faces out of~$\sC_i$.) The marked corner of the added light face is set arbitrarily.

\begin{figure}[ht]
		\psfrag{2}[][]{\textcolor{vert}{$e_2$}}
		\psfrag{C}[][]{\textcolor{red}{$\sC_2$}}
		\psfrag{-}[][]{\textcolor{red}{$\bwp_2^-$}}
		\psfrag{+}[][]{\textcolor{red}{$\bwp_2^+$}}
	\centering\includegraphics[width=.95\linewidth]{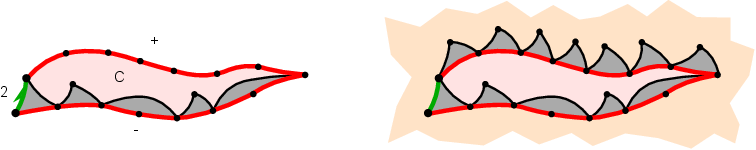}
	\caption{The connected component~$\sC_2$ and the corresponding $p$-hypermap (the irrelevant marked corner of the added light face is not depicted). The fact that the unbounded face is flawed implies the existence of a flawed face in the red area.}
	\label{pfirreg2}
\end{figure}

Finally, the added light face is necessarily a flawed face since, modulo~$p$, its degree is equal to $[\bwp_i^-]+1+(p-1)[\bwp_i^+] \equiv \vd(v, e_i^-)+1-\vd(v, e_i^+) \not\equiv 0$, by definition of~$e_i$ being irregular. As a result, there should be another flawed face in this $p$-hypermap, which is necessarily inside~$\sC_i$. The result follows.
\end{proof}

\section{Bijective interpretation}\label{secbij}

\subsection{Slit slide sew}\label{secsss}

\begin{figure}[ht!]
		\psfrag{5}[B][B][.8]{\unli{\textit{1. Simple path~\textcolor{red}{$\bwp$} linking~\textcolor{vert}{$c$} to~\textcolor{blue}{$c'$}}}}
		\psfrag{6}[B][B][.8]{\unli{\textit{2. Slit}}}
		\psfrag{7}[B][B][.8]{\unli{\textit{3. Slide}}}
		\psfrag{8}[B][B][.8]{\unli{\textit{4. Sew}}}
		
		\psfrag{f}[][][.8]{$f$}
		\psfrag{g}[][][.8]{$f'$}		
		\psfrag{c}[][][.7]{\textcolor{vert}{$c$}}
		\psfrag{d}[][][.7]{\textcolor{blue}{$c'$}}
		
		\psfrag{1}[][][.7]{\textcolor{red}{$e_1$}}
		\psfrag{2}[][][.7]{\textcolor{red}{$e_k$}}
		\psfrag{3}[][][.7]{\textcolor{red}{$\ell_1$}}
		\psfrag{4}[][][.7]{\textcolor{red}{$d_k$}}
		\psfrag{9}[][][.7]{\textcolor{red}{$\ell_2\sew d_1$}}
		
		\psfrag{y}[][][.7]{\textcolor{red}{$\bwp$}}
		\psfrag{l}[][][.8]{\textcolor{red}{$\bl$}}
		\psfrag{r}[][][.8]{\textcolor{red}{$\bd$}}
	\centering\includegraphics[width=.95\linewidth]{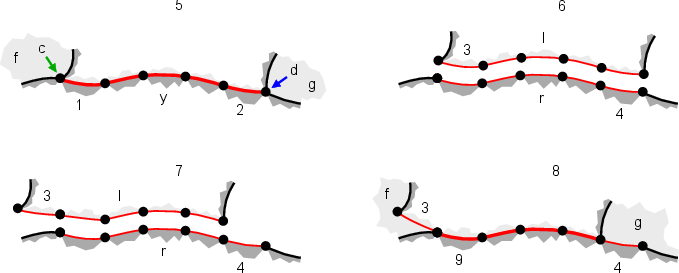}
	\caption{The slit-slide-sew operation on a $p$-hypermap. Note that, here, $\ell_1$ is not a dangling edge since~$c$ is not the corner preceeding~$e_1$.}
	\label{sss}
\end{figure}

Let us first describe the general operation at the heart of our construction. See Figure~\ref{sss}. Assume that, on some $p$-hypermap~$\m$, we have a simple path~$\bwp=(e_1,e_2,\dots,e_k)$ linking some corner~$c$ in some light face~$f$ to a different corner~$c'$ in some light face~$f'$ (which may possibly be equal to~$f$), that is, such that $e_1^-=c^+$ and $e_k^+=\smash{c'}^+$. We may then follow~$\bwp$, \textbf{entering from the corner~$c$} and \textbf{exiting through the corner~$c'$}. This creates a simple path on the sphere, starting inside the face~$f$ and finishing inside~$f'$. We may slit the sphere along this path, thus doubling the sides of the path. In the hypermap~$\m$, this doubles the path~$\bwp$, making up two copies, one incident to light faces and the slit, and one incident to dark faces and the slit. We denote by $\bl=(\ell_1,\dots,\ell_k)$ the former and by $\bd=(d_1,\dots,d_k)$ the latter.

Note that the data of~$\bwp$ is not sufficient to properly define this operation; one needs to know from which corner to enter~$\bwp$ in order to decide if an edge incident to~$e_1^-$ becomes incident whether to~$\ell_1^-$ or to~$d_1^-$. Similarly, one needs to know through which corner to exit~$\bwp$. 

We then sew back~$\bl$ onto~$\bd$ but only after sliding by one unit, in the sense that we match~$\ell_{i+1}$ with~$d_{i}$, for every $1\le i \le k-1$. For further reference, we denote by $\ell_{i+1}\sew d_{i}$ the resulting edge. Observe that, except from~$f$ and~$f'$, the faces are not altered by the process. Observe also that~$\ell_1$ and~$d_k$ are not matched with anything:
\begin{itemize}
	\item $d_k$ is still incident to the original dark face and is now also incident to~$f'$;
	\item $\ell_1$ is still incident to the original light face and is now also incident to~$f$.
\end{itemize}
Consequently, the result is no longer a $p$-hypermap since~$\ell_1$ is incident to light faces from both sides. However, in the case where~$\ell_1$ is actually a dangling edge (an edge with one extremity of degree~$1$), removing it provides a $p$-hypermap. See for instance Figure~\ref{figbijtwo}: the green dangling edge is removed in the final (bottom-right) picture, and actually replaced with a marked corner~$\tilde c$ in the face~$\tilde f_1$. This happens if and only if~$c$ is the corner preceding~$e_1$\,; this will always be the case in the present work.

\subsection{Face of degree two or more}\label{secdegtwo}

We now present the bijective interpretation for the identity~\eqref{eqpm} of Proposition~\ref{propmp1}.

\begin{figure}[ht!]
		\psfrag{f}[][][.7]{$f_1$}
		\psfrag{g}[][][.7]{$f_2$}
		\psfrag{i}[][][.7]{$\tilde f_1$}
		\psfrag{h}[][][.7]{$\tilde f_2$}
		
		\psfrag{m}[r][r][.8]{$\m$}
		\psfrag{1}[r][r][.8]{\textcolor{stepcol}{\textbf{1.--2.}}}
		\psfrag{3}[r][r][.8]{\textcolor{stepcol}{\textbf{3.}}}
		\psfrag{n}[r][r][.8]{$\tilde\m$}
		
		\psfrag{c}[][][.7]{\textcolor{blue}{$c$}}
		\psfrag{e}[][][.7]{\textcolor{vert}{$e$}}
		\psfrag{a}[][][.7]{\textcolor{vert}{$c_e$}}
		\psfrag{b}[][][.7]{\textcolor{blue}{$\tilde e$}}
		\psfrag{d}[][][.7]{\textcolor{vert}{$\tilde c$}}
		\psfrag{z}[][][.7]{\textcolor{red}{$\bgamma$}}
	\centering\includegraphics[width=.99\linewidth]{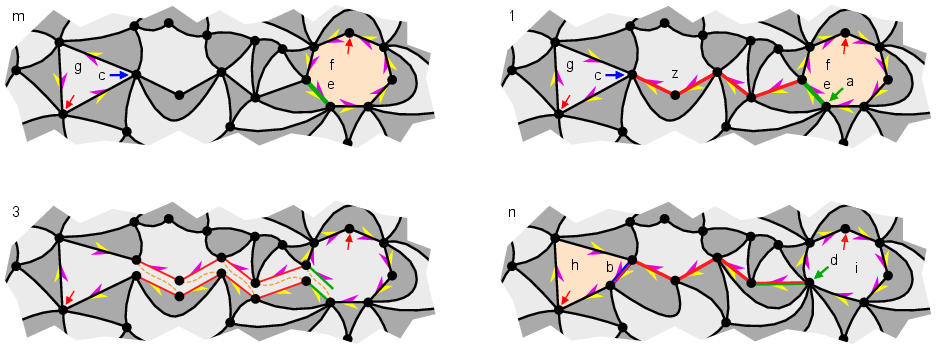}
	\caption{The involution $\Phi\colon\cH\to\cH$. Only the orientation around the faces of interest and along~$\bgamma$ are depicted. \textbf{Top left.} The input. \textbf{Top right.} We changed the orientation and defined the sliding path~$\bgamma$. \textbf{Bottom left.} We slit along the path. The dashed lines indicate to sew back after sliding. \textbf{Bottom right.} The output.}
	\label{figbijtwo}
\end{figure}

\paragraph{Involution.}
We start by defining a mapping~$\Phi$ on the set~$\cH$ of quadruples $(\co,\m,c,e)$, where
\begin{itemize}
	\item $\co$ is an orientation (either $\LL$ or $\LR$);
	\item $\m$ is a $p$-hypermap;
	\item $c$ is a distinguished corner of some light face;
	\item $e$ is a distinguished edge leaving~$c^{+}$ in the orientation~$\co$.
\end{itemize}
We break down the process into the following steps. See Figure~\ref{figbijtwo}.

\begin{description}[style=nextline]
\item[\textcolor{stepcol}{1. Reorientation}]
	From this point on, we convene to use the reverse orientation, which we denote by~$\tilde\co$.

\item[\textcolor{stepcol}{2. Sliding path}]
	We consider the corner~$c_{e}$ preceding~$e$ and the lightest geodesic~$\bgamma$ from~$c_{e}$ to~$c^+$.

\item[\textcolor{stepcol}{3. Slitting, sliding, sewing}]
	We slit, slide, sew along~$\bgamma$ from~$c_e$ to~$c$ as described in the previous section: along~$\bgamma$, the light side of an edge is now matched with the dark side of the previous edge. 

\item[\textcolor{stepcol}{4. Output}]
	The unmatched light side of the first edge of~$\bgamma$ yields a dangling edge; we remove it and denote the resulting corner by~$\tilde c$. We denote the edge corresponding to the unmatched dark side of the final edge of~$\bgamma$ by~$\tilde e$. We let~$\tilde\m$ be the resulting map. Finally, the output of the construction is the quadruple $\Phi(\co,\m,c,e)\de(\tilde\co,\tilde\m,\tilde c,\tilde e)$.
\end{description}

\begin{thm}\label{thminvol}
The mapping $\Phi\colon \cH \to \cH$ is an involution.
\end{thm}

\begin{proof}
We consider an element $(\co,\m,c,e)\in \cH$ and define $(\tilde\co,\tilde\m,\tilde c,\tilde e)\de\Phi(\co,\m,c,e)$. From the discussion in Section~\ref{secsss} and the definition of~$c_e$, we obtain that~$\tilde\m$ is a $p$-hypermap. In order to see that $(\tilde\co,\tilde\m,\tilde c,\tilde e)\in\cH$, it thus remains to see that~$\tilde e$ is leaving~$\tilde c^{+}$ in the orientation~$\tilde\co$.

Similarly to Section~\ref{secsss}, we denote the light side of the geodesic~$\bgamma$ by $\bl=(\ell_1,\dots,\ell_k)$ and its dark side by $\bd=(d_1,\dots,d_k)$. Let also~$c_{\tilde e}$ be the corner in~$\tilde\m$ that precedes~$\tilde e$ \textbf{in the orientation~$\co$}. Everything follows from the following claim.

\begin{claim}
In~$\tilde\m$, in the orientation~$\co$, the lightest geodesic from~$c_{\tilde e}$ to~$\tilde c^+$ is the path $(\tilde e, d_{k-1}\sew \ell_k, \dots, d_1 \sew \ell_2)$.
\end{claim}

Indeed, this implies that~$\tilde e$ is leaving~$\tilde c^{+}$ in the orientation~$\tilde\co$, so that $(\tilde\co,\tilde\m,\tilde c,\tilde e)\in\cH$. Furthermore, applying~$\Phi$ to $(\tilde\co,\tilde\m,\tilde c,\tilde e)$ means reverting~$\tilde\co$ back to~$\co$ and sliding back along the previous sliding path. More precisely, let us set $\tilde\bgamma=(\tilde e, d_{k-1}\sew \ell_k, \dots, d_1 \sew \ell_2)$ and denote by $\tilde\bl=(\tilde\ell_1,\dots,\tilde\ell_k)$ and $\tilde\bd=(\tilde d_1,\dots,\tilde d_k)$ its light and dark sides. For $2\le i \le k$, we thus have $\tilde\ell_{k+2-i}=\ell_{i}$ and $\tilde d_{k+1-i}=d_i$, so that we match back $\tilde\ell_{k+2-i}\sew \tilde d_{k+1-i}=\ell_i\sew d_i$. Furthermore, $\tilde d_k=d_1$ is unmatched, which amounts to matching it with the extra dangling edge~$\ell_1$ that we suppressed. Finally, the distinguished elements correspond.

\begin{figure}[ht]	
		\psfrag{m}[r][r][.8]{$\m$}
		\psfrag{n}[r][r][.8]{$\tilde\m$}
	
		\psfrag{e}[][][.7]{\textcolor{vert}{$e$}}
		\psfrag{a}[][][.7]{\textcolor{vert}{$c_e$}}
		\psfrag{d}[][][.7]{\textcolor{vert}{$\tilde c$}}
		\psfrag{b}[][][.7]{\textcolor{blue}{$\tilde e$}}
		\psfrag{h}[][][.7]{\textcolor{blue}{$c_{\tilde e}$}}
		\psfrag{c}[][][.7]{\textcolor{blue}{$c$}}
		\psfrag{y}[][][.7]{\textcolor{red}{$\tilde\bgamma$}}
		\psfrag{z}[][][.7]{\textcolor{red}{$\bgamma$}}
	
		\psfrag{l}[][][.8]{\textcolor{red}{$\bl$}}
		\psfrag{r}[][][.8]{\textcolor{red}{$\bd$}}
	
		\psfrag{w}[][][.7]{\textcolor{violet}{$\tilde\bwp$}}
		\psfrag{u}[][][.6]{\textcolor{violet}{$\boldsymbol{\le}$}}
		\psfrag{v}[][][.6]{\textcolor{violet}{$\bm{<}$}}
	\centering\includegraphics[width=.99\linewidth]{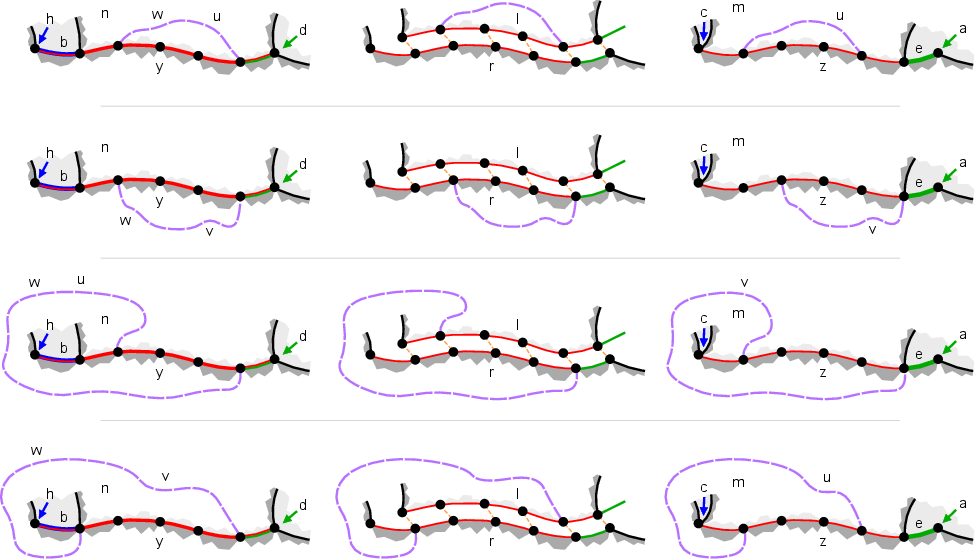}
	\caption{Proof of the fact that the sliding path~$\bgamma$ in~$\m$ used for applying~$\Phi$ becomes the (reverse of the) sliding path in the image~$\tilde\m$ used for applying~$\Phi$ a second time. \textbf{Left column.} If $\tilde\bwp\ne\tilde\bgamma$, it means that~$\tilde\bwp$ uses one of the purple dashed circumventions (each of the four lines corresponds to a possibility); if it starts on the light side, then it is shorter (symbol~\raisebox{.2ex}{\scalebox{.8}{$\boldsymbol\le$}}); if it starts on the dark side, then it is strictly shorter (symbol~\raisebox{.2ex}{\scalebox{.8}{$\boldsymbol<$}}). \textbf{Middle and right column.} Tracking it back to the original map~$\m$ yields a circumvention providing a better option than~$\bgamma$.}
	\label{pfinvol}
\end{figure}

Let us then show the claim, using Figure~\ref{pfinvol} as support. We argue by contradiction and assume that the lightest geodesic~$\tilde\bwp$ from~$\tilde e$ to~$\tilde c^+$ differs from~$\tilde\bgamma$. Then~$\tilde\bwp$ has to leave~$\tilde\bgamma$ at some point (either from its light side or from its dark side) and merge back at some other point (also either from its light side or from its dark side). This circumvention is either strictly shorter (a shortcut) or has the same length but comes before in the direction incoming edge, light face around the vertex where the splitting occurs. Tracing back this circumvention in~$\m$ shows that this would have provided a better option than~$\bgamma$, creating either a shortcut or a better circumvention, thus contradicting the definition of~$\bgamma$.
\end{proof}

\paragraph{Specialization.}
We now see how~$\Phi$ specializes into a bijection interpreting~\eqref{eqpm}. We let 
\[
\ba=(a_1,\ldots, a_r)\qquad\text{ and }\qquad\tilde \ba=(\tilde a_1,\ldots,\tilde a_r)\de( a_1-1, a_2+1, a_3,\ldots, a_r)
\]
be as in the statement of Proposition~\ref{propmp1}. Note that this means that $p$-hypermaps of type~$\ba$ are either $p$-constellations or quasi-$p$-constellations whose \textbf{first} face is flawed. Similarly, $p$-hypermaps of type~$\tilde\ba$ are either $p$-constellations or quasi-$p$-constellations whose \textbf{second} face is flawed.

We fix an orientation~$\co$ and define the following sets, whose cardinalities are respectively the left-hand side and the right-hand side of~\eqref{eqpm}, by Corollary~\ref{propq} (recall also the convention at the begining of Section~\ref{secprel} for distinguishing corners).

\medskip
\noindent\begin{minipage}[t]{.48\linewidth}
\begin{itemize}
\item We let~$\cM$ be the set of $p$-hypermaps of type~$\ba$ carrying 
\begin{itemize}
	\item one distinguished corner~$c$ in the \textbf{second} face,
	\item one distinguished edge~$e$ incident to the \textbf{first} face and leaving~$c^{+}$, for the \textbf{orientation}~$\co$.
\end{itemize}
\end{itemize}
\end{minipage}\hfill\begin{minipage}[t]{.5\linewidth}
\begin{itemize}
\item We let~$\tilde\cM$ be the set of $p$-hypermaps of type~$\tilde\ba$ carrying
\begin{itemize}
	\item one distinguished corner~$\tilde c$ in the \textbf{first} face,
	\item one distinguished edge~$\tilde e$ incident to the \textbf{second} face and leaving~$\smash{\tilde c}^{+}$, for the \textbf{reverse orientation}~$\tilde\co$.
\end{itemize}
\end{itemize}
\end{minipage}

\begin{flushright}
		\psfrag{f}[][]{$f_1$}
		\psfrag{g}[][]{$f_2$}
		\psfrag{h}[][]{$\tilde f_1$}
		\psfrag{i}[][]{$\tilde f_2$}
		\psfrag{1}[][][.9]{\textcolor{red}{$+1$}}
		\psfrag{9}[][][.9]{\textcolor{red}{$-1$}}
		\psfrag{M}[][]{$\tilde\cM$}
		\psfrag{N}[][]{$\cM$}
		\psfrag{c}[][]{\textcolor{blue}{$c$}}
		\psfrag{e}[][]{\textcolor{vert}{$e$}}
		\psfrag{b}[][]{\textcolor{blue}{$\tilde e$}}
		\psfrag{d}[c][t]{\textcolor{vert}{$\tilde c$}}
	\includegraphics[width=.9\linewidth]{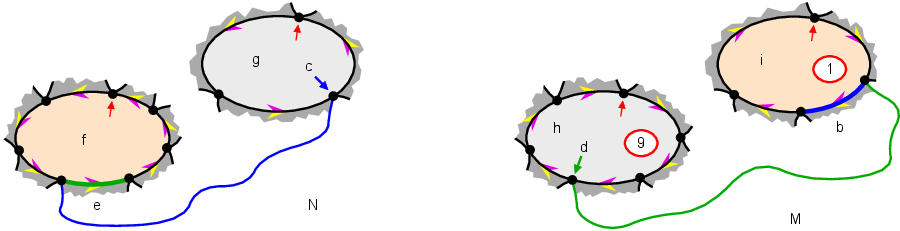}
\end{flushright}
\begin{center}
	\emph{Here, $p=3$, we are in the case~\ref{congruii} of Proposition~\ref{propmp1}, and $\co=\LL$.}
\end{center}

The pictograph above summarizes the definitions of~$\cM$ and~$\tilde\cM$. The red arrowheads denote the marked corners of the depicted faces; the half-arrowheads along the edges symbolize the convened orientation. Finally, the red~$\pm 1$ on the right shows the increase or decrease of the degree of the face in~$\tilde\cM$ in comparison with the one of the corresponding face in~$\cM$. In order to avoid confusion, we denote the first and second faces of maps in~$\cM$ by~$f_1$ and~$f_2$ as before, and use~$\tilde f_1$ and~$\tilde f_2$ instead, for maps in~$\tilde\cM$. The paths symbolize the fact that the edges are leaving the corners.

\begin{rem} Note that the convention on the orientation of edges is not the same in the definitions of the sets~$\cM$ and~$\tilde\cM$. This clearly bears no effects from an enumeration point of view but is of crucial importance for our bijections.
\end{rem}

\begin{corol}\label{cortwo}
The mapping~$\Phi$ specializes into a bijection from the set $\{(\co,\m,c,e)\,:\, (\m,c,e)\in\cM\}$ onto the set $\{(\tilde\co,\tilde\m,\tilde c,\tilde e)\,:\, (\tilde\m,\tilde c,\tilde e)\in\tilde\cM\}$, thus providing a bijection between~$\cM$ and~$\tilde\cM$.
\end{corol}

\subsection{Face of degree one}\label{secdegone}

We proceed to the bijective interpretation for the identity~\eqref{eqpm0} of Proposition~\ref{propface}, which works in a similar fashion as before.

\paragraph{Setting.}
Let $\ba=(1,a_2,\ldots, a_r)$ and $\tilde \ba=(\tilde a_2,\ldots,\tilde a_r)\de(a_2+1, a_3,\ldots, a_r)$ be tuples of positive integers, both with at most two coordinates not lying in~$p\N$. In order not to be confused by the index shift in~$\tilde\ba_2$, we denote the faces of $p$-hypermaps of type~$\tilde\ba$ by~$\tilde f_2$, \dots, $\tilde f_r$. In particular, $p$-hypermaps of type~$\tilde\ba$ are either $p$-constellations, or are quasi-$p$-constellations whose face~$\tilde f_2$ (the one with degree~$\tilde a_2$) is flawed. We fix an orientation~$\co$ and define the following sets, whose cardinalities are the sides of~\eqref{eqpm0}, again by Corollary~\ref{propq} for the right-hand side. 

\medskip
\noindent\begin{minipage}[t]{.46\linewidth}
\begin{itemize}
\item We let~$\cN$ be the set of $p$-hypermaps of type~$\ba$ carrying 
\begin{itemize}
	\item one distinguished corner~$c$ in the face~$f_2$.
\end{itemize}
\end{itemize}
\end{minipage}\hfill\begin{minipage}[t]{.52\linewidth}
\begin{itemize}
\item We let~$\tilde\cN$ be the set of $p$-hypermaps of type~$\tilde\ba$ carrying
\begin{itemize}
	\item one distinguished vertex~$\tilde v$,
	\item one distinguished edge~$\tilde e$ incident to~$\tilde f_2$ and leaving~$\tilde v$ for the orientation~$\tilde\co$.
\end{itemize}
\end{itemize}
\end{minipage}

\begin{center}
		\psfrag{f}[][]{$f_1$}
		\psfrag{g}[][]{$f_2$}
		\psfrag{i}[][]{$\tilde f_2$}
		\psfrag{1}[][][.9]{\textcolor{red}{$+1$}}
		\psfrag{M}[][]{$\tilde\cN$}
		\psfrag{N}[][]{$\cN$}
		\psfrag{c}[][]{\textcolor{blue}{$c$}}
		\psfrag{b}[][]{\textcolor{blue}{$\tilde e$}}
		\psfrag{v}[][]{\textcolor{vert}{$\tilde v$}}
	\includegraphics[width=.9\linewidth]{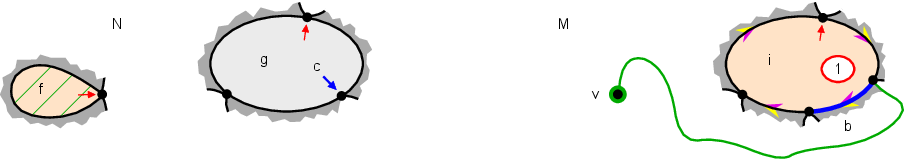} 
\end{center}

\begin{figure}[ht!]
		\psfrag{f}[][][.7]{$f_1$}
		\psfrag{g}[][][.7]{$f_2$}
		\psfrag{h}[][][.7]{$\tilde f_2$}
		
		\psfrag{m}[r][r][.8]{$\m$}
		\psfrag{1}[r][r][.8]{\textcolor{stepcol}{\textbf{1.--2.}}}
		\psfrag{3}[r][r][.8]{\textcolor{stepcol}{\textbf{3.}}}
		\psfrag{n}[r][r][.8]{$\tilde\m$}
		
		\psfrag{c}[][][.7]{\textcolor{blue}{$c$}}
		\psfrag{b}[][][.7]{\textcolor{blue}{$\tilde e$}}
		\psfrag{d}[][][.7]{\textcolor{vert}{$d_0$}}
		\psfrag{v}[][][.7]{\textcolor{vert}{$\tilde v$}}
		\psfrag{z}[][][.7]{\textcolor{red}{$\bgamma$}}
	\centering\includegraphics[width=.99\linewidth]{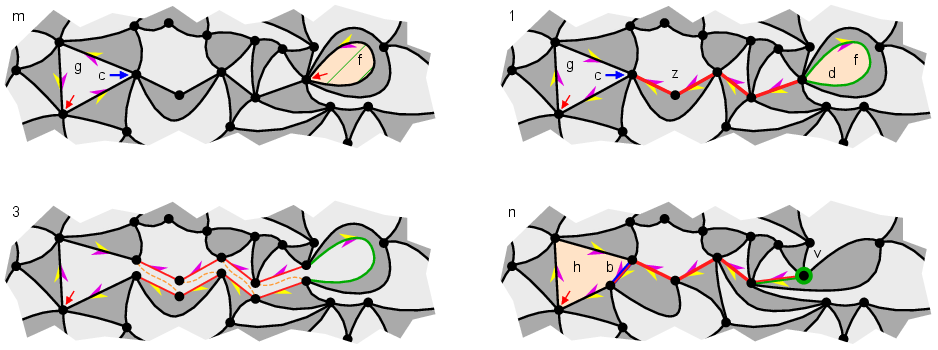}
	\caption{The bijection in the case of a degree~$1$-face, from~$\cN$ to~$\tilde\cN$.}
	\label{figbijone}
\end{figure}

We put~$f_1$ on the pictograph since we think of it as the ``missing'' distinguished element for~$\cN$. Note that we do not need to specify an orientation for maps in~$\cN$\,; we will however use the orientation~$\co$ for these maps in due time. The bijections between~$\cN$ and~$\tilde\cN$ can be thought of as degenerate versions of the one of the previous section. Here, we do not have an involution; we need to describe both mappings. We break them down into similar steps as above. See Figure~\ref{figbijone}.

\paragraph{Suppressing a face.}
We consider $(\m,c)\in \cN$.

\begin{description}[style=nextline]
\item[\textcolor{stepcol}{1. Reorientation}]
	From this point on, we use the reverse orientation~$\tilde\co$.

\item[\textcolor{stepcol}{2. Sliding path}]
	We forget the mark of the unique corner of~$f_1$ and consider the lightest geodesic~$\bgamma$ from this corner to~$c^+$.

\item[\textcolor{stepcol}{3. Slitting, sliding, sewing}]
	We denote by~$d_0$ the unique edge incident to~$f_1$. We slit, slide, sew along~$\bgamma$ from the unique corner of~$f_1$ to~$c$ as described in Section~\ref{secsss}, while furthermore matching the unmatched light side of the first edge with~$d_0$, that is, $\ell_1\sew d_0$ in the notation of Section~\ref{secsss}.

\item[\textcolor{stepcol}{4. Output}]
	We set $\Psi_{\boldsymbol-}(\m,c)\de(\tilde\m,\tilde v,\tilde e)$, where~$\tilde\m$ is the resulting map, $\tilde e$ is the edge corresponding to the unmatched dark side of the final edge of~$\bgamma$, and~$\tilde v$ is the origin of~$\bgamma$, that is, $\tilde v\de(\ell_1\sew d_0)^-$.
\end{description}

Note that the face~$f_1$ has been suppressed in the process.

\paragraph{Adding a face.}
We consider $(\tilde\m,\tilde v,\tilde e) \in \tilde\cN$.

\begin{description}[style=nextline]
\item[\textcolor{stepcol}{1. Reorientation}]
	We reverse the orientation, which means that we use the original orientation~$\co$.

\item[\textcolor{stepcol}{2. Sliding path}]
	We consider the corner~$c_{\tilde e}$ preceding~$\tilde e$ and the lightest geodesic~$\tilde\bgamma$ from~$c_{\tilde e}$ to~$\tilde v$.

\item[\textcolor{stepcol}{3. Slitting, sliding, sewing}]
We slit~$\tilde\m$ along~$\tilde\bgamma$, entering from~$c_{\tilde e}$ and stopping at~$\tilde v$, \textbf{without disconnecting the map} at~$\tilde v$. Denoting as before by $\tilde\bl=(\tilde\ell_1,\dots,\tilde\ell_k)$ and $\tilde\bd=(\tilde d_1,\dots,\tilde d_k)$ the light and dark sides of~$\tilde\bgamma$, this means that $\tilde\ell_k^+=\tilde d_k^+$.

We then sew back~$\bl$ onto~$\bd$ as before, by matching~$\ell_{i+1}$ with~$d_{i}$, for every $1\le i \le k-1$. In this matching, $\tilde\ell_k^+$ is matched with $\tilde d_{k}^-$, so that~$\tilde d_k$ creates a loop, whose unmatched side encloses an extra face, which we denote by~$f_1$ and mark at its unique corner.

\item[\textcolor{stepcol}{4. Output}]
As usual, we replace the unmatched~$\tilde\ell_1$ with a corner, which we denote by~$c$. We let~$\m$ be the resulting map, and set $\Psi_{\boldsymbol+}(\tilde\m,\tilde v,\tilde e)\de(\m,c)$.
\end{description}

\begin{thm}
The mappings $\Psi_{\boldsymbol-}\colon \cN \to \tilde\cN$ and $\Psi_{\boldsymbol+}\colon \tilde\cN \to \cN$ are well defined and inverse bijections.
\end{thm}

\begin{proof}
The proof closely follows that of Theorem~\ref{thminvol} and is left to the reader.
\end{proof}

\section{Uniform sampling}\label{secsample}

\paragraph{Algorithm.}
Recall the definition of the sequence~$\ba$ of types defined by~\eqref{eqseq} in the introduction, as well as Figure~\ref{figrow}. Recall also that we may easily sample a uniform $p$-constellation~$\m^0$ of type $\ba^0=(np)$. Let us then explain how to sample from there a sequence of $p$-hypermaps~$\m^0$, \dots, $\m^{np-a_r}$ such that, for every $1\le j\le np-a_r$, the $p$-hypermap~$\m^j$ is uniformly distributed among $p$-hypermaps of type~$\ba^j$.
\begin{itemize}
	\item Sample~$\m^0$.	    
	\item
	Assuming~$\m^{j-1}$ has been sampled, we sample~$\m^{j}$ as follows. 
	\begin{itemize}
		\item We randomly pick an orientation~$\co$.
		\item In~$\m^{j-1}$, we distinguish either a corner or a vertex:
		\begin{itemize}
			\item[\ref{samplea}] if a coordinate of~$\ba^{j-1}$ is increased by~$+1$, we distinguish a corner chosen uniformly at random in the corresponding face of~$\m^{j-1}$;
			\item[\ref{sampleb}] otherwise, we distinguish a vertex chosen uniformly at random in~$\m^{j-1}$.
		\end{itemize}
		\item Then, we distinguish an edge chosen uniformly at random among those incident to the face corresponding to the decreasing coordinate of~$\ba^{j-1}$ that leave the distinguished corner or vertex.
		\item We apply the mapping of Section~\ref{secdegtwo} in case~\ref{samplea} or that of Section~\ref{secdegone} in case~\ref{sampleb} in order to sample a $p$-hypermap~$\m^{j}$ of type~$\ba^j$ from~$\m^{j-1}$, the distinguished corner or vertex, and the distinguished edge.
		\item In the resulting $p$-hypermap~$\m^{j}$, we forget the distinguished elements coming from the bijection.
	\end{itemize}
\end{itemize}

As the number of ways to select the desired distinguished elements in~$\m^{j-1}$ only depends on the type~$\ba^{j-1}$, distinguishing elements does not bias the uniform probability. And since bijections preserve the uniform distribution, the fact that~$\m^{j-1}$ is uniformly distributed among the $p$-hypermaps of type~$\ba^{j-1}$ implies that~$\m^{j}$ is uniformly distributed among the $p$-hypermaps of type~$\ba^{j}$. By induction, we obtain that, for every $1\le j\le np-a_r$, the $p$-hypermap~$\m^j$ is indeed uniformly distributed among $p$-hypermaps of type~$\ba^j$. In particular, $\m^{np-a_r}$ is a uniform $p$-constellation or a quasi-$p$-constellation of type~$\ba$, as desired.

\paragraph{Practical implementation.}
In order to implement the algorithm in practice, one needs to 
\begin{enumerate}
	\item find all the edges incident to a given face that leave a given vertex;
	\item find a lightest geodesic from a corner to a vertex.
\end{enumerate}
To do so, one may use a directed version of Dijkstra's algorithm. For the first task, we can compute all the the distances to the vertex and select the proper edges. For the second task, we may compute all the distances to a the target vertex. Then, finding the lightmost geodesic amounts to select the leftmost or rightmost edge at each step of the path among those bringing us closer to the target vertex.

From Dijkstra's algorithm analysis, it comes that both tasks may be done in $\bO(np\log(np))$ time. On top of this, updating the data of the map should require a time linear in the number of edges (so $\bO(np)$); therefore, going from~$\m^{j-1}$ to~$\m^{j}$ can be done in $\bO(np\log(np))$ time. This slit-slit-sew operation has to be done $np - a_r$ times in order to obtain a $p$-hypermap of type $\ba = (a_1,\ldots,a_r)$ from a $p$-hypermap of type $(np)$. As a result, we obtain an overall complexity in $\bO(n^2 p^2\log(np))$.


\bibliographystyle{alpha}
\bibliography{main}

\begin{thebibliography}{BCCGF24}

\bibitem[ARS97]{AlReSc97}
L.~Alonso, J.~L. R\'{e}my, and R.~Schott.
\newblock A linear-time algorithm for the generation of trees.
\newblock {\em Algorithmica}, 17(2):162--182, 1997.

\bibitem[BCCGF24]{BoChChGF24}
Valentin Bonzom, Guillaume Chapuy, S\'{e}verin Charbonnier, and Elba
  Garcia-Failde.
\newblock Topological recursion for {O}rlov-{S}cherbin tau functions, and
  constellations with internal faces.
\newblock {\em Comm. Math. Phys.}, 405(8):Paper No. 189, 53, 2024.

\bibitem[Bet14]{Bet14}
J{\'e}r{\'e}mie Bettinelli.
\newblock Increasing forests and quadrangulations via a bijective approach.
\newblock {\em J. Combin. Theory Ser. A}, 122(0):107--125, 2014.

\bibitem[Bet19]{Bet19short}
J{\'e}r{\'e}mie Bettinelli.
\newblock Slit-slide-sew bijections for bipartite and quasibipartite plane
  maps.
\newblock In {\em 31th {I}nternational {C}onference on {F}ormal {P}ower
  {S}eries and {A}lgebraic {C}ombinatorics ({FPSAC} 2019)}, volume 82B of {\em
  S\'em. Lothar. Combin.}, pages Art. 82, 1--12. SLC, 2019.

\bibitem[Bet20]{Bet20}
J{\'e}r{\'e}mie Bettinelli.
\newblock Slit-slide-sew bijections for bipartite and quasibipartite plane
  maps.
\newblock {\em Electron. J. Combin.}, 27(3):Paper 3.4, 23, 2020.

\bibitem[BK24]{BeKo24short}
J{\'e}r{\'e}mie Bettinelli and Dimitri Korkotashvili.
\newblock Slit-slide-sew bijections for constellations and quasiconstellations.
\newblock In {\em 36th {I}nternational {C}onference on {F}ormal {P}ower
  {S}eries and {A}lgebraic {C}ombinatorics ({FPSAC} 2024)}, volume 91B of {\em
  S\'em. Lothar. Combin.}, pages Art.90, 1--12. SLC, 2024.

\bibitem[BMS00]{BMSc00}
Mireille Bousquet-M\'{e}lou and Gilles Schaeffer.
\newblock Enumeration of planar constellations.
\newblock {\em Adv. in Appl. Math.}, 24(4):337--368, 2000.

\bibitem[CF14]{CoFu14}
Gwendal Collet and {\'E}ric Fusy.
\newblock A simple formula for the series of constellations and
  quasi-constellations with boundaries.
\newblock {\em Electron. J. Combin.}, 21(2):Paper 2.9, 27, 2014.

\bibitem[DPS14]{DuPoSc14}
Enrica Duchi, Dominique Poulalhon, and Gilles Schaeffer.
\newblock Bijections for simple and double hurwitz numbers.
\newblock {\em Preprint,
  \href{http://arxiv.org/abs/1410.6521}{\textup{\nolinkurl{arXiv:1410.6521}}}},
  2014.

\bibitem[Eze78]{Eze78}
Cloyd~L. Ezell.
\newblock Branch point structure of covering maps onto nonorientable surfaces.
\newblock {\em Trans. Amer. Math. Soc.}, 243:123--133, 1978.

\bibitem[GJ83]{GoJa83}
Ian~P. Goulden and David~M. Jackson.
\newblock {\em Combinatorial enumeration}.
\newblock Wiley-Interscience Series in Discrete Mathematics. John Wiley \&
  Sons, Inc., New York, 1983.
\newblock With a foreword by Gian-Carlo Rota, A Wiley-Interscience Publication.

\bibitem[Hur91]{Hur91}
Adolf Hurwitz.
\newblock Ueber {R}iemann'sche {F}l\"{a}chen mit gegebenen
  {V}erzweigungspunkten.
\newblock {\em Math. Ann.}, 39(1):1--60, 1891.

\bibitem[Tut62]{Tut62sli}
William~T. Tutte.
\newblock A census of slicings.
\newblock {\em Canad. J. Math.}, 14:708--722, 1962.

\end{thebibliography}

\end{document}